\documentclass[11pt]{article}
\usepackage{amssymb}

\def\hang{\hangindent\parindent}
\def\textindent#1{\indent\llap{#1\enspace}\ignorespaces}
\def\re{\par\hang\textindent}
\def\QED{\hfill{$\Box$}}
\def \r{\rightarrow}\def\T#1{\widetilde #1}\def\op{\oplus}

\def\mapdown#1{\llap{$\vcenter {\hbox {$\scriptstyle #1$}}$}
                                \Bigg\downarrow}
\def\mapright#1#2{\smash{\mathop{\longrightarrow}\limits^{#1}_{#2}}}

\def\mapup#1{\Bigg\uparrow\rlap{$\vcenter {\hbox  {$\scriptstyle #1$}}$}}
\def\v5{\vskip .5truecm}
\def\OV#1{\overline {#1}}
\def\B{{\cal B}} \def\F{{\cal F}}
\def\LM{{\bf LM}}  
\def\LH{{\bf LH}}
\def\LR{\langle X\rangle} 
\def\LRY{\langle 1-X\rangle}\def\LRT{\langle 1-t\rangle}
\def\FRAC#1#2{\displaystyle{\frac{#1}{#2}}} 
\def\SUM^#1_#2{\displaystyle{\sum^{#1}_{#2}}}

\title{The General PBW Property\thanks{Project supported by 
the National Natural Science Foundation of China (10571038).}}
\author{Huishi Li\\
{\small Department of Applied Mathematics}\\
{\small College of Information Science and Technology}\\
{\small Hainan University}\\
{\small Haikou, 570228, China}}
\date{}

\begin{document}
\maketitle
\begin{center}
\begin{minipage}{120mm} 
{\small {\bf Abstract.}
For ungraded quotients of an arbitrary $\mathbb{Z}$-graded ring, we define the general 
PBW property, that covers the classical PBW property and the $N$-type PBW property 
studied via the $N$-Koszulity by several authors ([BG1], BG2], [FV]). 
In view of the noncommutative Gr\"obner basis theory, we conclude that every ungraded quotient of 
a path algebra (or a free algebra) has the general PBW property. We remark that  
an earlier result of Golod [Gol] concerning Gr\"obner bases can be used to give a homological characterization 
of the general PBW property in terms of Shafarevich complex. Examples of application 
are given.  }
\end{minipage}\end{center}{\parindent=0pt\par
}
\footnote{2000 Mathematics Classification. Primary: 16W70; Secondary: 16Z05.}
\vskip 1truecm
\section*{0. Introduction} 
Let $K\LR$ be the free associative algebra on a set of noncommuting variables $X$ over a 
field $K$, and let $K\LR =\oplus_{p\in\mathbb{N}}K\LR_p$ be the decomposition of $K\LR$ by its 
homogeneous components $K\LR_p$ spanned by words of length $p\ge 0$. Then $K\LR$ has the natural filtration 
$FK\LR=\{ F_pK\LR\}_{p\in\mathbb{N}}$ with $F_pK\LR =\oplus_{i\le p}K\LR_i$. For a $K$-subspace $P\subset 
F_NK\LR$, $N\ge 2$, let $\langle P\rangle$ be the two-sided ideal of 
$K\LR$ generated by $P$ and write $A=K\LR /\langle P\rangle$. Then $FK\LR$ induces a filtration 
$FA=\{ F_pA\}_{p\in\mathbb{N}}$ on $A$, where $F_pA=(F_pK\LR +\langle P\rangle )/\langle P\rangle$, 
that defines the associated graded $K$-algebra $G(A)=\oplus_{p\in\mathbb{N}}G(A)_p$ with $G(A)_p=F_pA/F_{p-1}A$. 
Since 
$$\FRAC{K\LR_p\oplus F_{p-1}K\LR}{(F_pK\LR\cap\langle P\rangle)+F_{p-1}K\LR}= 
\FRAC{F_pK\LR}{(F_pK\LR\cap\langle P\rangle)+F_{p-1}K\LR}~\mapright{\cong}{}~G(A)_p,\quad p\in\mathbb{N},$$
there is the natural graded epimorphism $\phi$: $K\LR\r G(A)$. On the other hand, consider 
$P_N=\{ f\in P~|~f\not\in F_{N-1}K\LR\}\subseteq P$. Then every $f\in P_N$ has a unique 
presentation $f=f_N+f_{N-1}+\cdots +f_{N-s}\in P$ with $f_N\in K\LR_N$, $f_{N-j}\in K\LR_{N-j}$ 
and $f_N\ne 0$. Write $\LH (f)=f_N$ for each $f\in P_N$. Then $\LH (f)\in (\langle P\rangle\cap F_NK\LR )+F_{N-1}K\LR$. 
Thus, if $\langle\LH (P_N)\rangle$ denotes the graded two-sided ideal of $K\LR$ 
generated by $\LH (P_N)=\left\{\LH (f)~\Big |~f\in P_N\right\}$, then $\langle\LH (P_N)\rangle\subseteq$  
Ker$\phi$. It follows that the canonical graded epimorphism $\pi$: $K\LR\r \OV A=K\LR /\langle\LH (P_N)\rangle$ 
yields naturally a graded epimorphism $\rho$: $\OV A\r G(A)$ such that the following diagram commutes
$$\begin{array}{ccc} K\LR&\mapright{\pi}{}&\OV A\\
\mapdown{\phi}&\swarrow\scriptstyle{\rho}&\\ G(A)&&\end{array}$$
Actually, the property that {\it $\rho$ is an isomorphism} is an analogue of the classical 
PBW (abbreviation of Poincar\`e-Birkhoff-Witt) theorem for enveloping algebras of Lie algebras.  
For $N=2$, Braverman and Gaitsgory [BG1] studied this isomorphism problem posed by Joseph Bernstein. Applying graded 
deformations to both graded hochschild cohomology and Koszul algebras, they obtained a 
PBW theorem as follows.\v5{\parindent=0pt
{\bf Theorem A} ([BG1] Theorem 0.5) Suppose that $P$ satisfies\par
(I) $P\cap F_1K\LR =\{ 0\}$ and \par
(J) $(F_1K\LR\cdot P\cdot F_1K\LR )\cap F_2K\LR =P$.\par
If the quadratic algebra $\OV A=K\LR /\langle\LH (P_2)\rangle$ is Koszul in the classical sense, then 
$\rho$ is an isomorphism. }\v5
If we call the PBW property studied in [BG1] the 2-{\it type PBW property} for the reason that 
$P\subset F_2K\LR$, then generally for $N\ge 2$,  the $N$-{\it type PBW property} 
was studied in the very recent work [FV] and [BG2] respectively. 
Gunnar Floystad and Jon Eivind Vatne dealt with the $N$-type PBW property in [FV] for 
deformations of $N$-Koszul $K$-algebras and the obtained $N$-type 
PBW theorem states that
{\parindent=0pt\v5
{\bf Theorem B} ([FV] Theorem 4.1) Suppose that the graded algebra $\OV A=K\LR /\langle\LH (P_N)\rangle$ is an $N$-Koszul algebra in 
the sense of [Ber]. Then $\rho$ is an isomorphism if and only if 
$$\langle P\rangle\cap F_NK\LR =P.$$
While Roland Berger and Victor Ginzburg dealt with the $N$-type PBW property in [BG2] 
for ungraded quotients of the tensor algebra over a 
Von Neumann regular ring $K$ and an $N$-type PBW theorem was obtained as well.
\v5
{\bf Theorem C} ([BG2] Theorem 3.4)  Suppose that $P$ satisfies\par
(a) $P\cap F_{N-1}K\LR =\{ 0\}$ and \par
(b) $(P\cdot K\LR_1+K\LR_1\cdot P )\cap F_NK\LR =P$.\par
If the graded left $K$-module Tor$^A_3(K_A,{~}_AK)$ is concentrated in degree $N+1$, 
then $\rho$ is an isomorphism.
\v5
As a consequence, an extension of the $N$-Koszulity [Ber] to nonhomogeneous algebras was  
realized through the $N$-type PBW property in [BG2]. }
{\parindent=0pt\v5
{\bf Remark} Note that we have stated Theorems A--C in the language of the present paper. For instance, 
in [FV], the algebra $\OV A$ is a given $N$-Koszul algebra defined by homogeneous elements of degree $N$, 
and the algebra $A$ is a deformation of $\OV A$ 
such that its associated graded algebra is exactly $\OV A$. }
\v5
In this paper, for ungraded quotients of an arbitrary  $\mathbb{Z}$-graded ring, 
we define first the general PBW property, that covers the classical PBW property and 
the $N$-type PBW property studied via the $N$-Koszulity 
in the literature. This is reached in section 1 after a clear picture of ungraded quotients 
vs graded quotients is established (Theorem 1.6). 
In section 2, we focus on ungraded quotients of path algebras (including free algebras) and realize 
the general PBW property by means of Gr\"obner bases. We remark in section 3 that an earlier result of Golod [Gol] concerning Gr\"obner bases 
can be used to give a homological characterization of the general PBW property (for positively graded 
algebra) in terms of Shafarevich complex. Finally in section 4, some examples of applications of sections 1--2 are discussed. \par
Here we point out that the main idea and principal method used in section 1 
and section 2 were announced in Chapter III of [Li], where similar results were discussed only for 
quotients of finitely generated free algebras but the general PBW property for ungraded quotients of 
graded algebras was not exposed.\par
Throughout this paper, by a graded ring we mean an associative $\mathbb{Z}$-graded ring with unity 1. Let  
$B=\oplus_{p\in\mathbb{Z}}B_p$ be a graded ring. If $B_i=0$ for all $i<0$, then we say that $B$ is positively graded 
and write $B=\oplus_{p\in\mathbb{N}}B_p$.  
We adopt the conventional notion on graded rings and call an element $F_p\in B_p$ a homogeneous 
element of degree $p$. Thus, if $f=f_p+f_{p-1}+\cdots +f_{p-s}$ with $f_p\in B_p$, $F_{p-j}\in B_{p-j}$ and $f_p\ne 0$,  
then we say that $f$ has degree $p$ and write $d(f)=p$. Let $I$ be an ideal of $B$. Then 
$I$ is a graded ideal if and only if $I=\oplus_{p\in\mathbb{Z}}(B_p\cap I)$ if and only if 
the quotient ring $B/I=\oplus_{p\in\mathbb{Z}}\left (B_p+I/I\right )$. Unless otherwise stated, 
all graded ring (module) homomorphisms are of degree 0.
\v5
\section*{1. The General PBW Property}
In this section we introduce and characterize the general PBW property for ungraded quotients of an arbitrary 
$\mathbb{Z}$-graded ring. Since such a property is defined for filtered rings, from both 
a structural and a computational viewpoint (see Proposition 1.7 and Theorem 2.2), it is natural to 
bring both the associated graded ring and the Rees ring into the data considered. To begin with, 
let us review some necessary results on filtered rings and their associated graded objects. \par
Let $A$ be a $\mathbb{Z}$-filtered associative ring with filtration $FA$: 
$$\cdots\subset F_{p-1}A\subset F_pA\subset F_{p+1}A\subset\cdots, \quad p\in\mathbb{Z},$$
where each $F_pA$ is an abelian subgroup of $A$ (if $A$ is a $K$-algebra over some commutative ring $K$, 
then $F_pA$ is a $K$-submodule of $A$) such that $1\in F_0A$, $A=\cup_{p\in\mathbb{Z}}F_pA$, and 
$F_pAF_qA\subset F_{p+q}A$ for all $p,q\in\mathbb{Z}$. $FA$ induces two graded structures, that is, 
the associated graded ring $G(A)$ of $A$ which is defined as  $G(A)=\oplus_{p\in\mathbb{Z}}G(A)_p$ with $G(A)_p=F_pA/F_{p-1}A$, and the 
Rees ring $\T A$ of $A$ which is defined as $\T A=\oplus_{p\in\mathbb{Z}}F_pA$. Write $X$ for 
the homogeneous element of degree 1 in $\T A_1=F_1A$ represented by 1, which is usually called 
the {\it canonical element} of $\T A$. Then $X$ 
is contained in the center of $\T A$ and is not a divisor of 0. Consider 
the ideal $\langle 1-X\rangle=(1-X)\T A$, respectively $\langle X\rangle =X\T A$, of $\T A$ generated by 
$1-X$, respectively by $X$. Then it is well-known that 
$$\T A/\langle 1-X\rangle\cong A,\quad\quad \T A/X\T A\cong G(A).\leqno{(*)}$$\par
On the other hand, Let $B=\op_{n\in\mathbb{Z}}B_n$ 
be a graded ring and $X$ a homogeneous element of degree 1, that is, $X\in B_1$. Suppose
that $X$ is contained in the center of $B$ and is not a divisor of 0.  Put $\Lambda=B/\langle 1-X\rangle$, 
where $\langle 1-X\rangle$ is the ideal of $B$ generated by $1-X$.  
Note that if $b_p\in B_p$ is a homogeneous element of degree $p$, then $Xb_p\in
B_{p+1}$ because $X$ is of degree 1. Thus
 $b_p=Xb_p +(1-X)b_p$ implies $B_p+\langle 1-X\rangle\subset B_{p+1}+\langle
1-X\rangle$. Consequently, the $\mathbb{Z}$-gradation on $B$ induces naturally a $\mathbb{Z}$-filtration
$F\Lambda$ on $\Lambda$: 
$$F_n\Lambda =\displaystyle{\frac{B_p+\LRY}{\LRY}} ,\quad p\in\mathbb{Z}.$$
\v5{\parindent=0pt
{\bf 1.1. Proposition} (see [LVO]) With notation as above, the following statements hold.\par
(i) $\LRY$ does not contain any nonzero homogeneous element of $B$.\par
(ii) The associated graded ring $G(\Lambda )$ of $\Lambda$ with respect to $F\Lambda$ 
is isomorphic to $B/XB$ under graded ring homomorphism.\par
(iii) The Rees ring $\T{\Lambda}$ of $\Lambda$ with respect to $F\Lambda$ is 
isomorphic to $B$ under graded ring homomorphism. In particular, $X$ corresponds to    
the canonical element of $\T{\Lambda}$.}\par\QED
\vskip 6pt
For the remainder of this section, let $R=\op_{n\in\mathbb{Z}}R_n$ be an arbitrary graded ring.\vskip 6pt
Consider the polynomial ring $R[t]$ over $R$ in one
commuting variable $t$. Then the onto ring homomorphism $\phi$: $R[t]\r R $ defined by 
$\phi (t)=1$ has Ker$\phi =\LRT$, the ideal of $R[t]$ generated by $1-t$. Hence $R\cong R[t]/\LRT$.
Since $R[t]$ has the  mixed gradation, that is, $R[t]=\op_{p\in\mathbb{Z}}R[t]_p$ with
$$R[t]_p=\left\{\left. \displaystyle{\sum_{i+j=p}}F_it^j~\right |~F_i\in R_i,
~j\ge 0
 \right \} ,\quad p\in \mathbb{Z} ,$$
for each $f\in R$, there exists a homogeneous element $F\in
R[t]_p$, for some $p$, such that $\phi (F)=f$. More precisely, if $f=f_p+f_{p-1}+
\cdots +f_{p-s}$ where $f_p\in R_p$,  $f_{p-j}\in R_{p-j}$ and $f_p\ne 0$, 
then $f^*=f_p+tf_{p-1}+\cdots +t^sf_{p-s}$ is a homogeneous element in $R[t]_p$ satisfying $\phi (f^*)=f$.
{\parindent=0pt
\v5
{\bf 1.2. Definition} (i) For any $F\in R[t]$, write $F_*=\phi
(F)$. $F_*$ is called the {\it dehomogenization} of $F$ with respect to
$t$.\par
(ii) For an element $f\in R$, if $f=f_p+f_{p-1}+\cdots +f_{p-s}$ with 
$f_p\in R_p$, $f_{p-j}\in R_{p-j}$ and $f_p\ne 0$, then the
homogeneous element $f^*=f_p+tf_{p-1}+\cdots +t^sf_{p-s}$ in 
$R[t]_p$ is called the {\it homogenization} of $f$ with respect to
$t$.\par
(iii) If $I$ is a two-sided ideal of $R$, then we let $\langle I^*\rangle$ stand
for the graded two-sided ideal of $R[t]$ generated by $I^*=\{ f^*~|~f\in
I\}$. $\langle I^*\rangle$ is called the {\it homogenization ideal} of $I$ with respect to $t$.}
\v5
With definition and notation as above, the following 1.3--1.4 may be found in [Li].
{\parindent=0pt\v5
{\bf 1.3. Lemma} (i) For $F,G\in R[t]$, $(F+G)_*=F_*+G_*$,
$(FG)_*=F_*G_*$.\par
(ii) For $f,g\in R$, $(fg)^*=f^*g^*$, $t^s(f+g)^*=t^rf^*+t^hg^*$, where
$r=d(g)$, $h=d(f)$, and $s=r+h-d(f+g)$.\par
(iii) For any $f\in R$, $(f^*)_*=f$.\par
(iv) If $F$ is a homogeneous element of degree $p$ in $R[t]$, and if
$(F_*)^*$ is of degree $q$, then $t^r(F_*)^*=F$, where $r=p-q$.\par
(v) If $I$ is a two-sided ideal of $R$, then each homogeneous element 
$F\in \langle I^*\rangle$ is of the form $t^rf^*$ for some $r\in\mathbb{N}$ and $f\in I$.\par\QED
\v5
{\bf 1.4. Proposition} Let $I$ be a proper two-sided ideal of $R$.
Then the map
$$\begin{array}{cccc}
\alpha :&R[t]/\langle I^*\rangle&\mapright{}{}&R/I\\
&F+\langle I^*\rangle&\mapsto&F_*+I,
\end{array} \quad\quad F\in R[t]$$
is an onto ring homomorphism  with
Ker$\alpha =\langle 1-\overline{t}\rangle$, where $\overline{t}$ denotes the coset of $t$
in $R[t]/\langle I^*\rangle$. Moreover, $\overline{t}$ is not a divisor of 0  in $R[t]/\langle I^*\rangle$, 
and hence $\langle 1-\overline{t}\rangle$
does not contain any nonzero homogeneous element of $R[t]/\langle I^*\rangle$ .}\par\QED
\v5
Consider the natural grading filtration $FR$ on $R$ which is defined by the 
the abelian subgroups
$$F_pR =\op_{i\le p}R,\quad p\in\mathbb{Z}.$$
Let $I$ be a proper two-sided ideal of $R$ and $A=R/I$. Then $FR$
induces the quotient filtration $FA$ on $A$:
$$F_pA=(F_pR +I)/I,\quad p\in\mathbb{Z},$$
that defines two graded structures: the associated graded ring $G(A)=\op_{p\in\mathbb{Z}}G(A)_p$ 
of $A$ with $G(A)_p=F_pA/F_{p-1}A$,  and the Rees ring $\T A=\op_{p\in\mathbb{Z}}\T A_p$ of $A$ with 
$\T A_p=F_pA$. The proposition below shows that $G(A)$ and $\T A$ may be 
 determined by $\langle I^*\rangle$. 
{\parindent=0pt\v5
{\bf 1.5. Proposition} With notation as before, there are graded
ring isomorphisms:\par
(i) $\T A\cong R[t]/\langle I^*\rangle$, and\par
(ii) $G(A)\cong R[t]/(\langle t\rangle +\langle I^*\rangle)$, where $\langle
t\rangle$ denotes the ideal of $R[t]$ generated by $t$.
\vskip 6pt
{\bf Proof} Put $B=R[t]/\langle I^*\rangle$, $B_p=(R[t]_p+\langle I^*\rangle )/\langle I^*\rangle 
=\overline{R[t]_p}$, $p\in\mathbb{Z} $. Then  $\overline{t}$ is a homogeneous element of degree 1 
in $B$, and by Proposition 1.4, it is not a divisor of 0. Hence $B$ is isomorphic to 
the Rees ring $\T{\Lambda}$ of the filtered ring $\Lambda =B/\langle 1-\OV t\rangle$, where 
$$F_p\Lambda =\FRAC{B_p+\langle 1-\OV t\rangle}{\langle 1-\OV t\rangle}=
\frac{\overline{R[t]_p}+\langle 1-\OV{t}\rangle}{\langle 1-\OV{t}\rangle},\quad p\in\mathbb{Z},$$
and moreover,  $B/\OV tB\cong G(\Lambda )$. On the other hand, it is not difficult to see that 
the ring homomorphism $\alpha$: $R[t]/\langle I^*\rangle\r
R/I=A$ defined in Proposition 1.4 yields isomorphisms of abelian groups:
$$F_p\Lambda =\frac{\overline{R[t]_p}+\langle 1-\OV{t}\rangle}{\langle 1-\OV{t}\rangle}
~\mapright{}{}~\displaystyle{\frac{\op_{i\le p}R_i+I}{I}} =F_pA, \quad
p\in\mathbb{Z},$$
which extend to define a graded ring isomorphism 
$$\T{\alpha}:\quad B\cong \T{\Lambda}=\bigoplus_{p\in\mathbb{Z}}
\frac{\overline{R[t]_p}+\langle 1-\OV{t}\rangle}{\langle 1-\OV{t}\rangle}
\mapright{}{} \bigoplus_{p\in\mathbb{Z}}F_pA=\T A.$$
But note that under $\alpha$ we have $t+\langle I^*\rangle\mapsto 1+I$. Thus, under 
the graded ring isomorphism $\T{\alpha}$ we have $\OV t\mapsto X$, the canonical element of $\T A$. 
It follows from the formula $(*)$ given in the beginning of this section and 
Proposition 1.1 that (i) and (ii) hold.\QED}
\v5
Further, we present $G(A)$ as a graded quotient of $R$ by finding its defining ideal 
clearly. To this end, for $f\in R$ we denote by $\LH (f)$ the {\it leading 
homogeneous part} of $f$, that is, if $f=f_p+f_{p-1}+\cdots +f_{p-s}$ 
with $f_p\in R_p$,  $f_{p-j}\in R_{p-j}$ and $f_p\ne 0$, 
then $\LH(f)=f_p$. Thus, if $S$ is a subset of $R$, then we put 
$$\LH (S)=\left\{ \LH (f)~\Big |~f\in S\right\}$$
and write $\langle\LH (S)\rangle$ for the graded two-sided ideal generated by $\LH (S)$ in $R$. \par
Since $G(A)=\oplus_{p\in\mathbb{Z}}G(A)_p$, where $G(A)_p=F_pA/F_{p-1}A$ with $F_pA=(F_pR+I)/I$,  
there are canonical ismorphisms of abelian groups
$$\FRAC{R_p\oplus F_{p-1}R}{(I\cap F_pR)+F_{p-1}R}=\FRAC{F_pR}{(I\cap F_pR)+F_{p-1}R}
~\mapright{\cong}{}~G(A)_p,\quad p\in\mathbb{Z} .\leqno{(1)}$$
It follows that the natural epimorphisms of abelian groups 
$$\phi_p:\quad R_p~\mapright{}{}~\FRAC{R_p\oplus F_{p-1}R}{(I\cap F_pR)+F_{p-1}R},\quad p\in\mathbb{Z},$$
extend to define a graded epimorphism
$$\phi :\quad R~\mapright{}{}~G(A).$$
{\parindent=0pt\v5
{\bf 1.6. Theorem} With the convention made above, we have Ker$\phi =\langle\LH (I)\rangle$, and hence 
$G(A)\cong R/\langle\LH (I)\rangle$. 
\vskip 6pt
{\bf Proof} It is sufficient to prove the equalities
$$\hbox{Ker}\phi_p=\langle\LH (I)\rangle\cap R_p,\quad p\in\mathbb{Z}.$$
Suppose $f_p\in$ Ker$\phi_p$. Then $f_p\in (I\cap F_pR)+F_{p-1}R$. If $f_p\ne 0$, then 
as $f_p$ is a homogeneous element of degree $p$, we have $f_p=\LH (f)$ for some $f\in I\cap F_pR$. 
This shows that $f_p\in\langle\LH (I)\rangle\cap R_p$. Hence Ker$\phi_p\subseteq\langle\LH (I)\rangle\cap R_p$. 
Conversely, suppose $f_p\in\langle\LH (I)\rangle\cap R_p$. Then $f_p=\sum g_i\LH (f_i)h_i$, 
where $g_i$ and $h_i$  are homogeneous elements. Let $f_i=\LH (f_i)+f_i'$ where 
deg$(f_i')<$ deg$(f_i)$. Then $f_p=\sum g_if_ih_i-\sum g_if_i'h_i$ with 
$\sum g_if_ih_i\in I\cap F_pR$ and $\sum g_if_i'h_i\in F_{p-1}R$. This shows that $f_p\in 
(I\cap F_pR)+F_{p-1}R$, that is, $f_p\in$ Ker$\phi_p$. Hence, $(\langle\LH (I)\rangle\cap R_p)\subseteq$ Ker$\phi_p$. 
Summing up, we conclude the desired equalities. \QED}
\v5
Now, let $\F$ be an {\it arbitrary} subset  of the ideal $I$. Then by the foregoing discussion, 
$$\langle\LH(\F )\rangle\subseteq\langle\LH (I)\rangle =\hbox{Ker}\phi .$$
It follows that the canonical graded epimorphism $\pi$: $R\r\OV A=R/\langle\LH (\F )\rangle$ 
yields naturally a graded epimorphism $\rho$: $\OV A\r G(A)$ such that the following diagram 
commutes
$$\begin{array}{ccc} R&\mapright{\pi}{}&\OV A\\
\mapdown{\phi}&\swarrow\scriptstyle{\rho}&\\ G(A)&&\end{array}$$
If we set $R=K\LR$, $I=\langle P\rangle$, and $\F =P_N$ as in section 0,  
then the property that {\it $\rho$ is an isomorphism} gives exactly the $N$-{\it type 
PBW Property}. Instead of giving our definition of the general PBW property immediately 
by using the phrase ``$\rho$ is an isomorphism", let us see first how  
Theorem 1.6 reveals the essential feature of this property. 
{\parindent=0pt\v5
{\bf 1.7. Proposition} Let $\F$ be an arbitrary subset of the ideal $I$ and $\OV A=R/\langle\LH (\F )\rangle$. 
With the convention made above, the following statements are equivalent.\par
(i) The natural graded epimorphism $\rho$: $\OV A\r G(A)$ is an isomorphism.\par
(ii)  $\langle\LH (I)\rangle=\langle\LH (\F )\rangle$.\par
(iii) $\F$ is a set of generators for $I$ that has the property: every $f\in I$ has a presentation $f=\sum g_jf_jh_j$, where $g_j,h_j\in 
R$ and $f_j\in\F$, such that $d(g_j)+d(f_j)+d(h_j)\le ~\hbox{deg}(f)$ for all $g_jf_jh_j\ne 0$.\par
(iv) $\langle I^*\rangle =\langle\F^*\rangle$ and hence $\T A\cong R[t]/\langle\F^*\rangle$, where $\F^*=\{ f^*~|~f\in \F\}$. \vskip 6pt
{\bf Proof} (i) $\Leftrightarrow$ (ii) By the construction of $\rho$, this equivalence is clear.\par 
(ii) $\Leftrightarrow$ (iii) Suppose $\langle \LH (\F )\rangle=\langle \LH (I)\rangle$. 
If $f\in I$ with $d(f)=p$, then since $\LH (f)$ is a homogeneous element we have $\LH (f)=\sum g_j\LH (f_j)h_j$ 
for some homogeneous elements $g_j, h_j\in R$, $f_j\in\F$, and
$d(g_j)+d(\LH (f_j))+d(h_j)=d(g_j)+d(f_j)+d(h_j)=p$. Now the element $f'=f-\sum g_jf_jh_j\in I$ has
$d(f')<p$, so we may repeat the same procedure for $f'$. Since $d(f)=p$ is finite, 
after a finite number of reduction steps we obtain a presentation $f=\sum g_jf_jh_j$ 
where $g_j,h_j$ are homogeneous elements of $R$, $f_j\in \F$ and $d(g_j)+d(f_j)+d(h_j)\le p$ for
 all $i$. It follows that (iii) holds.}\par
Conversely, suppose (iii) holds. Then it is easy to see that 
for any $f\in I$, $\LH (f)=\sum\LH (g_j)\LH (f_j)\LH (h_j)$ for some $g_j,h_j\in R$,
$f_j\in\F$. Hence $\langle \LH (\F )\rangle =\langle \LH (I)\rangle$. {\parindent=0pt\par
(iv) $\Leftrightarrow$ (iii) To prove this equivalence, first recall and bear in mind 
that if $f_p+f_{p-1}+\cdots +f_{p-s}=f\in R$ 
with $f_p\in R_p$, $f_{p-j}\in R_{p-j}$ and $f_p\ne 0$, then $f^*=f_p+tf_{p-1}+\cdots 
+t^sf_{p-s}$. Consequently, $d(f)=d(f^*)=p$ and $\LH (f)=\LH(f^*)=f_p$.}\par
Suppose (iv) holds. Then for $f\in I$ with $d(f)=p$, we have $f^*\in\langle I^*\rangle$. Hence, 
$f^*=\sum G_jf^*_jH_j$, in which $f_j^*\in\F^*$, $G_j$ and $H_j$ are homogeneous elements of $R[t]$ and $d(G_j)+d(f^*_j)+d(H_j)=p$ 
whenever $G_jf^*_jH_j\ne 0$. It follows from Lemma 1.3 that
$$f=(f^*)_*=\sum G_{j*}(f_j^*)_*H_{j*}=\sum G_{j*}f_jH_{j*},$$
where $d(G_{j*})+d(f_j)+d(H_{j*})\le p$ whenever $G_{j*}f_jH_{j*}\ne 0$. This shows that (iii) holds.\par
Conversely, suppose (iii) holds. To reach (iv), we need only to consider homogeneous elements. 
If $F\in\langle I^*\rangle$ is a homogeneous element, then by Lemma 1.3, $F=t^rf^*$ for some integer $r\ge 0$ and 
some $f\in I$. Suppose $f=\sum_jh_jf_jg_j$. Then $d(h_j^*)+d(f_j^*)+d(g_j^*)\le d(f^*)$. We may use 
Lemma 1.3 and the assumption (iii) to start a reduction procedure as follows.{\parindent=0pt \par 
{\bf Begin}
$$\begin{array}{l} f^*-\sum_jh_j^*f_j^*g_j^*=t^{r_1}m_1^*+t^{r_2}m_2^*+\cdots~\hbox{with}
\\
r_j>0,~m_j\in I,~\hbox{and}~d(t^{r_j}m_j^*)\le d(f^*)~\hbox{ for~
all}~m_j.\end{array}$$
For each $m_j^*\in I^*$,  where
 $m_j=\sum_ih_{i_j}f_{i_j}g_{i_j}$, since $d(h_{i_j}^*)+d(f_{i_j}^*)+d(g_{i_j}^*)\le d(m_j^*)$, 
so go to \par
{\bf Next}
$$\begin{array}{l} m_j^*-\sum_ih_{i_j}^*f_{i_j}^*g_{i_j}^*=t^{r_{1_j}}m_{1_j}^*+t^{r_{2_j}}m_{2_j}^*+\cdots~\hbox{with}\\
r_{k_j}>0,~m_{k_j}\in I,~\hbox{and}~d(t^{r_{k_j}}m_{k_j}^*)\le
d(m_j^*)~\hbox{ for~ all}~m_{k_j}.\end{array}$$
As $d(f^*)$ is finite, after a finite number of steps we may 
reach  $f^*\in\langle\F^*\rangle$, in particular,
$f^*=\sum_jh_j^*f_j^*g_j^*$ with $d(h_j^*)+d(f_j^*)+d(g_j^*)\le d(f^*)$
for all $j$. (Since the ideal $I$ considered should be a proper
ideal, the dehomogenization operation on $R[t]$ guarantees that the final result of the 
reduction procedure cannot be an
expression like $\sum_{\ell}t^{\ell}$.) This proves the conclusion of (iv).
 \QED}\v5
Proposition 1.7 tells us that if $\rho$ is an isomorphism, then the subset $\F$ is necessarily 
a set of generators for the ideal $I$, that is, $\F$ is not really ``arbitrary".
{\parindent=0pt\v5
{\bf 1.8. Definition} Let $R$, $I$ and $A=R/I$ be as before, and let $\F$ be a set of 
generators for the ideal $I$. The ring $A$ is said to have the {\it general PBW property}, 
if one of the equivalent conditions in Proposition 1.7 is satisfied.
\v5
{\bf Remark} (i) Clearly, Definition 1.8 covers the $N$-type 
PBW property, and it is also obvious that if $I$ is a graded ideal, then this definition 
becomes trivial.  If $I$ is not a graded ideal, then,  just like verifying the 
sufficient conditions for the $N$-type PBW property in Theorems A--C of [BG1], [FV] and [BG2], 
any of the equivalent conditions in Proposition 1.7 is not easy to be verified.   
We will see in next section that for ideals of a path algebra (or a free algebra), 
Gr\"obner bases with respect to a certain gradation-preserving monomial ordering 
can realize Proposition 1.7 effectively.\par
(ii)  Suppose that $I$ is ungraded, or equivalently, $A=R/I$ is 
not a graded ring. Then, except reaching a unified definition for the PBW property, 
the importance of Theorem 1.6 may also be indicated from a viewpoint of lifting structures. 
For instance, if $R$ is a finitely generated free algebra
or a finitely generated path algebra, and  if the ring $R/\langle\LH (I)\rangle$ is one of the folowing type: 
a domain, a Noetherian ring, an Artinian ring, a graded semisimple ring, a ring with 
finite global dimension, an Auslander regular ring, a ring with classical standard PBW-basis, etc, then 
$A=R/I$ is a ring of the same type at the ungraded level, and moreover, all properties 
listed may be lifted to the Rees ring $\T A$ of $A$ (see [LVO]).}
\v5
\section*{2. Gr\"obner Basis Means The General PBW Property} \def\G{{\cal G}}
In this section, we realize the general PBW property for quotients of path algebras (including 
free algebras) by means of Gr\"obner bases.  In principle, 
as the {\it Noncommutative Buchberger Algorithm} ([Mor], [Gr]) produces a (finite or infinite) 
Gr\"obner basis for each two-sided ideal of a path algebra (or a free algebra), we may 
say, from both a theoretical and a practical viewpoint, that every ungraded quotient 
of a path algebra (or a free algebra) has the general PBW property. \vskip 6pt
Before starting the main text of this section, let us explain briefly why path algebra 
is our first choice. Let $K$ be a field and $Q$ a finite 
directed graph (or a quiver). Recall that the path algebra $KQ$ is defined to be the $K$-algebra with the $K$-basis the set 
of finite directed paths in $Q$, where the vertices of $Q$ are viewed as paths of length 0, and 
the multiplication in $KQ$ is induced by multiplication of paths. 
Note that the free associative $K$-algebra on $n$ noncommuting variables is isomorphic 
to the path algebra $KQ$ where $Q$ has one vertex and $n$ loops, and hence, every finitely generated $K$-algebra 
is of the form $KQ/I$, where  $I$ is a two-sided ideal of $KQ$. It is well-known from 
representation theory that every finite dimensional $K$-algebra is Morita equivalent 
to an algebra of the form $KQ/I$ if $K$ is algebraically 
closed; and since $KQ$ has the natural gradation defined by 
the lengths of the paths, quotients of path algebra over $K$ include graded $K$-algebras 
$A=\oplus_{i\ge 0}A_i$, where $A_0$ is a product of a finite number of copies of $K$, 
each $A_i$ is a finite dimensional $K$-vector space and $A$ is generated in degree 0 and 1, 
that is, for $i,j\ge 0$, $A_iA_j=A_{i+j}$. It is also known that an algebra is $N$-Koszul 
in the sense of [Ber] if and only if it is a quotient of a path algebra by an ideal generated 
by homogeneous elements of degree $N$ and its Yoneda algebra is generated in degree 0, 1 and 2. 
Thus, our choice of path algebra has a big 
generality. In particular, every path algebra $KQ$ holds a well-developed 
Gr\"obner basis theory. So, defining relations of a quotient of $KQ$ may be studied algorithmically,  
and this advantage enables us to reach the main result of this section. \par
For a general theory on noncommutative Gr\"obner bases, 
the reader is referred to, for example, [Mor], [Gr] and [Li]. 
\v5
To maintain the notation of section 1, let us write $R=KQ$ and  
use the natural positively graded structure $R=\oplus_{p\in\mathbb{N}}R_p$ on $R$, 
where the gradation is defined by the lengths of paths in $R$. 
Let $I$ be a two-sided ideal of $R$ and $A=R/I$. Then $A$ has the filtration $FA$ 
induced by the grading filtration $FR$ on $R$. Let $G(A)$ and $\T A$ be the 
associated graded algebra and the Rees algebra of $A$ defined by $FA$, respectively. 
Then by Proposition 1.5 and Theorem 1.6, there are isomorphisms of graded $K$-algebras: 
$G(A)\cong R/\langle\LH (I)\rangle$, $\T A\cong R[t]/\langle I^*\rangle$.\par
From now on in this section we fix an admissible system $(R ,\B ,\succeq_{gr})$,  
that is, $\B$ is the $K$-basis of $R$ consisting of monomials (finite directed 
paths), and $\succeq_{gr}$ is some {\it graded monomial ordering} on $\B$, for example, 
the graded lexicographic ordering. If $f\in R$, 
$f=\sum\lambda_iu_i$ where $\lambda_i\in K$ and $u_i\in\B$, then we write
$\LM(f)=\max\{ u_i~|~\lambda_i\ne 0\}$ for the {\it leading monomial} of $f$. For 
a subset $S$ of $R$ we put 
$$\LM (S)=\left\{\LM (f)~\Big |~f\in S\right\}  $$
and write $\langle\LM (S)\rangle$ for the two-sided monomial ideal of $R$ generated by $\LM (S)$. 
Recall that a subset $\G \subset I$ is called a {\it Gr\"obner basis} for the two-sided ideal $I$ if 
$$\langle\LM (I)\rangle =\langle\LM (\G )\rangle .$$
{\parindent=0pt\v5
{\bf 2.1. Theorem}  Let $\G$  be a Gr\"obner basis for the two-sided 
ideal $I$ in $R$ with respect to $\succeq_{gr}$, and $A=R/I$. With notation maintained from section 1, 
the following statements hold.\par
(i) $\langle\LH(I)\rangle =\langle\LH(\G )\rangle$, and hence $G(A)\cong R/\langle\LH(\G )\rangle$, 
that is, the algebra $A$ has the general PBW property in the sense of Definition 1.8. \par
(ii) $\langle I^*\rangle =\langle\G^*\rangle$, and hence $\T A\cong R[t]/\langle\G^*\rangle$.\vskip 6pt
{\bf Proof} Since $\G$ is a Gr\"obner basis for $I$, it is well-known that for $f\in I$, starting with 
$\LM (f)$, the equality $\langle\LM (I)\rangle = \langle\LM (\G )\rangle$ (or a division on $f$ by $\G$) 
yields inductively a Gr\"obner presentation 
$$f=\sum\lambda_ju_jg_jv_j,\quad \lambda_j\in K,~u_j,v_j\in\B,~g_j\in\G,$$
in which $u_j\LM (g_j)v_j\ne 0$ and $\LM (f) \succeq_{gr}\LM (u_jg_jv_j)$. But 
note that the monomial ordering $\succeq_{gr}$ preserves gradation. It follows that 
$\LM (f)$ comes from $\LH (f)$ and they have the same degree at the graded level. Consequently, 
in the  Gr\"obner presentation of $f$ obtained above we also have 
$$d(f)\ge d(u_j)+d(g_j)+d(v_j)~\hbox{for all}~u_jg_jv_j.$$
This shows that Proposition 1.7(iii) is  satisfied. Hence (i) and (ii) hold.\QED}
\v5
In computational algebra it is a well-known fact that, starting with a set of homogeneous 
elements, a homogeneous Gr\"obner basis may be 
obtained in a more effective way.  At this point, in addition to its own independent 
interest, the next theorem will be helpful 
in realizing the general PBW property algorithmically.\par 
Consider the $k$-basis $$\B (t)=\left\{ wt^r~\Big |~w\in \B ,~r\ge 0\right\}$$ 
for $R[t]$. Then the monomial ordering $\succeq_{gr}$ on $\B$ extends to a
monomial ordering on $\B (t)$, again denoted $\succeq_{gr}$, as follows:
$$w_1t^{r_1}\succ_{gr}w_2t^{r_2}~\hbox{if~and ~only~if~} w_1\succ_{gr}
w_2,~\hbox{or}~w_1=w_2~\hbox{and}~r_1>r_2.$$\par
With the definition made above, we have $X_j\succ_{gr} t^r$ for all $j\in\Lambda$
and all $r\ge 0$, and a Gr\"obner basis theory holds in $R[t]$ exactly as  in $R$. 
{\parindent=0pt\v5
{\bf 2.2. Theorem}  Let $I$ be a two-sided ideal of $R$ and $\G\subset I$. With notation as before, the following statements 
are equivalent:\par
(i) $\G$ is a Gr\"obner basis of $I$ in $R$;\par
(ii) $\LH (\G )$ is a Gr\"obner basis of $\langle \LH (I)\rangle$ in $R$.\par
(iii) $\G^*$ is a Gr\"obner basis of $\langle I^*\rangle$ in $R[t]$; \par
\vskip 6pt
{\bf Proof}  First note that $\succeq_{gr}$ and  the definition of homogenization 
yield the following equalities:
$$\left\{\begin{array}{l}  \LM (f)=\LM (\LH (f)),~f\in R, \\
\LM (f^*)=\LM (f),~f\in R.\\
\end{array}\right.\leqno{(**)}$$ 
(i) $\Leftrightarrow$ (ii) Every element of $\langle\LH (I)\rangle$ 
has a presentation of the form $\sum\lambda_iu_i\LH(f_i)v_i$, where $\lambda_i\in K$, $u_i,v_i\in \B$ and
$f_i\in I$. Consequently, by the formula $(**)$ above, the desired equivalence follows from 
the equivalence below: for $f\in I$, $g_j\in\G$, 
$$\LM (f)=u\LM (g_j)u\Leftrightarrow \LM(\LH (f))=u\LM (\LH (g_j))v.$$\par
(i) $\Leftrightarrow$ (iii) Suppose (i) holds. Noticing that  $\langle I^*\rangle$ is a graded ideal, it needs only 
to consider homogeneous elements of $\langle I^*\rangle$. Let $F\in \langle I^*\rangle$ be a nonzero homogeneous element. 
Then by Lemma 1.3, $F=t^rf^*$, where $r\ge 0$ and $f\in R$. It follows from the foregoing formula $(**)$ that 
$$\begin{array}{rcl} \LM (F)&=&t^r\LM (f^*)\\
&=&t^r\LM (f)\\
&=&t^ru\LM (g_j)v,~\hbox{for some}~u,v\in\B,~g_j\in\G ,\\
&=&t^ru\LM (g_j^*)v.\end{array}$$
This shows that $\langle\LM (\langle I^*\rangle )\rangle \subseteq \langle\LM (\G^* )\rangle$, and 
hence the equality holds. Therefore, $\G^*$ is a Gr\"obner basis of $\langle I^*\rangle$. }\par
Conversely, suppose (iii) holds. Then for any $f\in I$, by the formula $(*)$ we have
$$\begin{array}{rcl} \LM (f)&=&\LM (f^*)\\
&=&u\LM (g^*_j)v,~\hbox{for some}~u,v\in\B ,~g_j\in\G ,\\
&=&u\LM (g_j)v.\end{array}$$
This shows that $\langle\LM (I)\rangle\subseteq \langle\LM (\G)\rangle$, and hence 
the equality holds. Therefore, $\G$ is a Gr\"obner basis of $I$.\QED
\v5
\section*{3. A Characterization in Terms of Shafarevich Homology}
In this section, we remark that the general PBW property for {\it positively $\mathbb{N}$-graded algebra}  
can be characterized by the first homology of the Shafarevich complex. This is based on an earlier work of 
Golod [Gol] in which {\it standard bases} (including Gr\"obner bases in path algebras and 
free algebras) was studied by means of Shafarevich homology and 
the classical Koszulity was involved in the commutative case. 
To understand this, note first that for a positively $\mathbb{N}$-graded $K$-algebra $R=\oplus_{p\in\mathbb{N}}R_p$, 
where $K$ is a field, if we adopt the notion and notation of [Gol] by setting $\Gamma=\mathbb{N}$ 
and using the grading $\mathbb{N}$-filtration $FR$ as the $\Gamma$-filtration, then 
the property $\langle\LH (I)\rangle=\langle\LH (\F )\rangle$ (Proposition 1.7(iii)) is just an 
analogue of the definition for a {\it standard basis} (including Gr\"obner basis) $\F$ in $R$. 
It turns out that our general PBW property also has a homological 
characterization, as to which, we mention now as follows. \vskip 6pt
Let $X=\{ x_j\}_{j\in J}$, and let $U$ denote the free $K$-algebra $K\LR$ on the set $X$. 
By definition, the {\it Shafarevich complex} relative to $F$, denoted Sh$(X|\F ,~R)$, is a complex of 
$R$-$R$-bimodules
$$\hbox{Sh}_n(X|\F ,~R)=\underbrace{R\otimes U\otimes R\otimes\cdots\otimes R\otimes U}_{n~\hbox{copies of}~ R\otimes U}\otimes R,\quad n\ge 0,$$
(tensor product is defined over $K$) and differentials
$$d_n:\quad\hbox{Sh}_n(X|\F ,~R)~\mapright{}{}~\hbox{SH}_{n-1}(X|\F ,~R)$$
with \def\OT{\otimes} \def\SH{\hbox{Sh}}
$$\begin{array}{l} d_n(a_0\OT x_{j_1}\OT a_1\OT\cdots\OT x_{j_{i-1}}\OT a_{i-1}\OT x_{j_i}\OT a_i\OT x_{j_{i+1}}\OT\cdots \OT a_{n-1}\OT x_{j_n}\OT a_n)\\
\\
=\SUM^n_{i=1}(-1)^{i-1}a_0\OT x_{j_1}\OT a_1\OT\cdots\OT x_{j_{i-1}}\OT (a_{i-1}f_{j_i}a_i)\OT x_{j_{i+1}}\OT\cdots \OT a_{n-1}\OT x_{j_n}\OT a_n),\end{array}$$\par
Consider the grading filtration $FR$ on $R$ as before. Then $FR$ induces a $\mathbb{N}$-filtration $F\hbox{Sh}(X|F~,R)$ on 
$\SH (X|\F ,~R)$, where for each $p,~n\ge 0$, $F_p\SH_n(X|\F ,~R)$ is the $K$-subspace spanned 
by 
$$a_0\OT x_{j_1}\OT a_1\OT\cdots\OT a_{n-1}\OT x_{j_n}\OT a_n$$ 
in which $a_i$s are homogeneous elements such that 
$$\hbox{deg}(a_0)+ \hbox{deg}(f_{j_1}+\hbox{deg}(a_1)+\cdots +\hbox{deg}(a_{n-1})+\hbox{deg}(f_{j_n})+\hbox{deg}(a_n)\le p.$$
With respect to this filtered structure, $d_n$ is a filtered homomorphism of degree 0, and hence, 
$\SH (X|\F ,~R)$ becomes a $\mathbb{N}$-filtered complex. It follows that there are two associated $\mathbb{N}$-graded 
complexes $G(SH )(X|\F ,~R)$ and $\T{\SH}(X|\F ,~R)$, where 
$$\bigoplus_{p\in\mathbb{N}}\FRAC{F_p\SH_n(X|\F ,~R)}{F_{p-1}\SH_n(X|\F ,~R)}=G(\SH )_n(X|\F ,~R)~\mapright{G(d_n)}{}~
G(\SH )_{n-1}(X|\F ,~R)=\bigoplus_{p\in\mathbb{N}}\FRAC{F_p\SH_{n-1}(X|\F ,~R)}{F_{p-1}\SH_{n-1}(X|\F ,~R)} ,$$
and 
$$\bigoplus_{p\in\mathbb{N}}F_p\SH_n(X|\F ,~R)=\T{\SH}_n(X|\F ,~R)
~\mapright{\T{d_n}}{}~\T{\SH}_{n-1}(X|\F ,~R)=\bigoplus_{p\in\mathbb{N}}F_p\SH_{n-1}(X|\F ,~R) .$$\par
Put $\LH (\F )=\{ \LH (f_j)\}_{j\in J}$. Then there is the Shafarevich complex $\SH (X|\LH (\F ),~R)$ 
relative to $\LH (\F )$ with differentials $D_n$, $n\ge 0$. Now, the natural graded surjective morphism 
$$\varphi :\quad \SH (X|\LH (\F ),~R)~\mapright{}{}~G(\SH )(X|\F ,~R)$$
and the canonical graded morphism
$$\psi :\quad \T{\SH}(X|\F ,~R)~\mapright{}{}~G(\SH )(X|\F ,~R)$$
induce homomorphisms of corresponding homology modules $\varphi_*$ and $\psi_*$, respectively.
Let $E_*(\SH (X|\LH (\F ),~R))$ denote the graded $R$-$R$-submodule $\varphi_*^{-1}(\hbox{Im}\psi_*)$ 
of $H_*(\SH (X|\LH (\F ),~R))$. Homogeneous elements in $E_*(\SH (X|\LH (\F ),~R))$ and 
the cycles representing them are called {\it extendable classes} and {\it cycles} respectively. 
Focusing on the first homology and tracing along the diagram
$$\begin{array}{ccccccc} \r&\SH_2(X|\LH (\F ),~R)_p&\mapright{D_{2p}}{}&\SH_1(X|\LH (\F ),~R)_p&
\mapright{D_{1p}}{}&\SH_0(X|\LH (\F ),~R)_p&\r\\
&\mapdown{\varphi_{2p}}&&\mapdown{\varphi_{1p}}&&\mapdown{\varphi_{0p}}&\\
\r&G(\SH )_2(X|\F ,~R)_p&\mapright{G(d_{2})_p}{}&G(\SH )_1(X|\F ,~R)_p&\mapright{G(d_{1})_p}{}&G(\SH )_0(X|\F ,~R)&\r\\
&\mapup{\psi_{2p}}&&\mapup{\psi_{1p}}&&\mapup{\psi_{0p}}&\\
\r&\T{\SH}_2(X|\F ,~R)_p&\mapright{\T{d_{2p}}}{}&\T{\SH}_1(X|\F ,~R)_p&\mapright{\T{d_{1p}}}{}&\T{\SH}_0(X|\F ,~R)_p&\r\end{array}$$
a homological characterization of the general PBW property is obtained as a special case of 
([Gol], Theorem 1).{\parindent=0pt\v5
{\bf 3.1. Theorem} With notation as before, the following statements are equivalent.\par
(i) $\langle\LH (\F )\rangle =\langle\LH (I)\rangle$, that is, $G(A)\cong R/\langle\LH (\F )\rangle$.\par
(ii) $E_1(\SH (X|\LH (\F ),~R))=H_1(\SH (X|\LH (\F ),~R))$.\par
(iii) The $R$-$R$-bimodule $H_1(\SH (X|\LH (\F ),~R))$ is generated by extendable classes.}
\v5
\section*{ 4. Examples}
The obvious application of previous sections 1--2 may be seen from Remark (ii) of section 1. 
In consideration of Koszulity, we finish this paper with two examples. \vskip 6pt
Let $R=\oplus_{p\in\mathbb{N}}R_p$ be a path algebra defined by a finite directed graph (or let $R$ be a finitely 
generated free algebra) over a field $K$, where the positive gradation is defined by the lengths of paths. 
If $I$ is generated by homogeneous elements of degree 2, 
then $A=R/I$ is called a {\it quadratic algebra}. One of the themes in the study of quadratic algebra 
has been the Koszulity (the well-known fact is that if $A=R/I$ is Koszul in the classical sense, 
then $I$ is generated necessarily by homogeneous elements of degree 2). Applying 
noncommutative Gr\"obner basis theory to $R$, if, with respect to a fixed monomial 
ordering $\prec$ on the standard $K$-basis $\B$ of $R$, the reduced Gr\"obner basis of $I$ (it always exists) 
consists of quadratic homogeneous elements, then $A$ is Koszul (for instance, see [GH]). Combined with 
the $N$-type PBW property, the $N$-Koszulity in the sense of [Ber] is generalized to 
ungraded quotients ([BG2], Definition 3.9), that is, taking the grading filtration $FR$ on $R$ 
into account, for $P\subset F_NR$, $N\ge 2$, 
and $I=\langle P\rangle$, the algebra $A=R/I$ is said to be Koszul if the graded 
algebra $R/\langle\LH (P_N )\rangle$ is $N$-Koszul and if the $N$-type PBW property holds (see 
section 0 for the notation used here).
{\parindent=0pt\v5
{\bf 4.1. Example} Let $I$ be an ideal of $R$, and let $\G =\{ g_j\}_{j\in J}$  be
 a Gr\"obner basis for $I$ with respect to some graded monomial ordering $\succeq_{gr}$ 
on $\B$. Consider the grading filtration $FR$ on $R$ and the induced filtration $FA$ on 
the quotient algebra $A=R/I$. With notation maintained from previous sections, the following 
statements hold.\par
(i) $A$ has the general PBW property in the sense of Definition 1.8, that is, the associated graded algebra 
$G(A)$ is isomorphic to $R/\langle\LH (\G )\rangle$. Moreover, $\LH (\G )$ is a 
Gr\"obner basis for $\langle\LH (I)\rangle$.\par
(ii) If $\G\subset F_2R$ and $\LH (\G )\ne 0$, then $G(A)$ is Koszul in the classical sense; 
If $\G\subset F_NR$ for $N\ge 2$ such that $\LH (\G )\ne 0$, then $A$ is Koszul in the sense of 
[BG2] whenever $R/\langle\LH (\G )\rangle$ is $N$-Koszul in the sense of [Ber].\par
(iii) The Rees algebra $\T A$ of $A$ is isomorphic to $R[t]/\langle\G^*\rangle$. Moreover, 
$\G^*$ is a Gr\"obner basis for $I^*$. \par
(iv) In the case that $\G\subset F_2R$, $\T A$ is Koszul in the classical sense. 
\v5
{\bf 4.2. Example} Consider any quadric solvable polynomial algebra $A$ studied 
in [Li] (in particular, Examples (i)--(vi) constructed in section 2 of Chapter III). 
Then $A$ has the following properties.\par
(a) $A$ has the general PBW Property in the sense of Definition 1.8 (indeed they 
all have classical standard PBW bases).\par
(b) With respect to its natural filtration $FA$ (induced by the grading filtration of a free algebra), 
$G(A)$ is Koszul in the classical sense.\par
(c) $A$ is Koszul in the sense of [BG2].\par
(d) The Rees algebra $\T A$ of $A$ with respect to $FA$ is a classical Koszul algebra.
\v5
{\bf Final remark} The result of Theorem 1.6 can be generalized to consider quotient algebras of 
a $\Gamma$-graded algebra $R$ with the $\Gamma$-grading filtration $\F R$, where $\Gamma$ is an ordered semigroup with a total ordering $\prec$. 
For instance, let $R=K\langle x_1,...,x_n\rangle$ be a finitely generated free $K$-algebra over a 
field $K$, and let $\B$ be the standard $K$-basis of $R$ consisting of words of length $\ge 0$. If 
$\prec$ is a  monomial ordering on $\B$, then, noticing now $R=\oplus_{u\in\B}Ku$ is $\B$-graded, 
we may consider the grading $\B$-filtration $\F R$ of $R$. Let $I$ be a two-sided ideal of $R$ and $A=R/I$. 
Then $\F R$ induces a $\B$-filtration $\F A$ for $A$ that defines the associated $\B$-graded 
algebra $G^{\F}(A)$ of $A$. In a similar way, we can reach an analogue of Theorem 1.6, 
that is, $G^{\F}(A)$ is isomorphic to the {\it monomial algebra} $R/\langle\LM (I)\rangle$, 
where $\LM (I)$ is the set of all leading monomials of $I$. If furthermore 
$I$ is generated by a Gr\"obner basis $\G$, then $G^{\F}(A)\cong R/\langle\LM (\G )\rangle$. A detailed 
discussion on this result and its applications will be given in a forthcoming paper.}
\v5
\centerline{References}{\parindent=1.2truecm\par
\re{[Ber]} R. Berger, Koszulity for nonquadratic algebras, {\it J. Alg}., 239(2001), 705--734.
\re{[BG1]} A. Braverman and D. Gaitsgory, Poincar\`e-Birkhoff-Witt theorem for quadratic algebras of Koszul 
type, {\it J. Alg.}, 181(1996), 315-328.
\re{[BG2]} R. Berger and V. Ginzburg, Symplectic Reflection algebras and non-homogeneous $N$-Koszul property, 
eprint arXiv:math.RA/0506093.
\re{[FV]}  G. Floystad and J. E. Vatne, PBW-deformations of $N$-Koszul algebras, \par
eprint arXiv:math/0505570.
\re{[GH]} E. Green and R. Q. Huang, Projective resolutions of straightening closed algebras generated by 
minors, {\it Adv. Math.}, 110(1995), 314-333.
\re{[Gol]} E. S. Golod, Standard bases and homology, in: {\it Some Current Trends in Algebra}, 
(Varna, 1986), LNM, 1352, Springer-Verlag, 1988, 88-95.
\re{[Gr]} E. Green, Noncommutative Gr\"obner bases and projective resolutions, in: {\it Computational Methods for 
Representations of Groups and Algebras} (Essen, 1997), Progr. Math., 173, Birkhuser, Basel, 1999, 29--60.
\re{[Li]} H. Li, {\it Noncommutative Gr\"obner Bases and Filtered-Graded Transfer}, LNM, 1795, Springer-Verlag, 2002.
\re{[LVO]} H. Li and F. Van Oystaeyen, {\it Zariskian Filtrations}, Kluwer Academic Publishers, 1996.
\re{[Mor]} T. Mora, An introduction to commutative and noncommutative Gr\"obner Bases, 
{\it Theoretic}\par {\it Computer Science}, 134(1994), 131--173.}

\end{document}